\documentclass[11pt, twoside]{article}
\usepackage{latexsym}
\usepackage{cite}
\usepackage{amsmath}\def\rad{\mbox{rad}}
\usepackage{amssymb}\def\End{\mbox{End}}
\usepackage[all]{xy}
\usepackage{amsfonts}
\usepackage{verbatim}
\usepackage{amsthm}\newcommand{\supp}{\mathsf{Supp}\hspace{.01in}}
\usepackage{mathrsfs}
\usepackage{epsfig}
\usepackage{xy}
\usepackage{array}
\usepackage{stmaryrd}
\usepackage{graphicx,color}
\usepackage{xcolor}
\usepackage[colorlinks=true,linkcolor=blue,citecolor=blue]{hyperref}
\usepackage{tikz}
\usetikzlibrary{arrows,calc}
\usepackage{etex}
\usepackage{mathdots}
\usepackage{float}
\usepackage{graphics}
\usepackage{pdflscape}
\usepackage{extarrows}
\usepackage{anysize,hyperref}
\input xypic
\xyoption{all}
\usepackage{bm}
\newcommand{\ind}{\mathsf{ind}\hspace{.01in}}
\usepackage[perpage,symbol]{footmisc}
\usepackage{setspace}
\def\U{\mathcal{U}}
\def\P{\mathcal{P}}
\def\I{\mathcal{I}}

\def\C{\mathscr{C}}
\def\E{\mathbb{E}}

\def\s{\mathfrak{s}}
\def\id{\mathrm{id}}
\def\op{^\mathrm{op}}
\def\Ab{\mathsf{Ab}}
\def\del{\delta}
\def\dr{\ar@{->}[r]}

\def\Hom{\mbox{Hom}}

\newcommand{\CC}{{\bf{C}}^{n+2}_{\C}}
\newcommand{\mr}{\hbox{\boldmath$\cdot$}}
\newcommand{\ov}{\overset}
\newcommand{\lra}{\longrightarrow}
\newcommand{\co}{\colon}
\newcommand{\uas}{^{\ast}}            %%% ^*
\newcommand{\sas}{_{\ast}}
\newcommand{\Xd}{\langle X_{\bullet},\del\rangle}  %%% <X,¦Ä>
\newcommand{\ush}{^\sharp}           %%% ^sharp
\newcommand{\ssh}{_\sharp}
\SelectTips{cm}{10}

\usepackage{setspace}
\begin{document}
\baselineskip=15pt
\title{\Large{\bf On the existence of Auslander-Reiten $\bm{n}$-exangles in\\[1mm] $\bm{n}$-exangulated categories\footnotetext{Jiangsheng Hu was supported by the NSF of China (Grant Nos. 12171206 and 11771212) and the Natural Science Foundation of Jiangsu Province (Grant No. BK20211358). Panyue Zhou was supported by the National Natural Science Foundation of China (Grant No. 11901190) and the Scientific Research Fund of Hunan Provincial Education Department (Grant No. 19B239).} }}
\medskip
\author{Jian He, Jiangsheng Hu, Dongdong Zhang and Panyue Zhou}

\date{}

\maketitle
\def\blue{\color{blue}}
\def\red{\color{red}}

\newtheorem{theorem}{Theorem}[section]
\newtheorem{lemma}[theorem]{Lemma}
\newtheorem{corollary}[theorem]{Corollary}
\newtheorem{proposition}[theorem]{Proposition}
\newtheorem{conjecture}{Conjecture}
\theoremstyle{definition}
\newtheorem{definition}[theorem]{Definition}
\newtheorem{question}[theorem]{Question}
\newtheorem{remark}[theorem]{Remark}
\newtheorem{remark*}[]{Remark}
\newtheorem{example}[theorem]{Example}
\newtheorem{example*}[]{Example}
\newtheorem{condition}[theorem]{Condition}
\newtheorem{condition*}[]{Condition}
\newtheorem{construction}[theorem]{Construction}
\newtheorem{construction*}[]{Construction}

\newtheorem{assumption}[theorem]{Assumption}
\newtheorem{assumption*}[]{Assumption}

\baselineskip=17pt
\parindent=0.5cm

\begin{abstract}
\begin{spacing}{1.2}
Let $\C$ be an $n$-exangulated category.
In this note, we show that if $\C$ is locally finite, then $\C$ has Auslander-Reiten $n$-exangles. This unifies
and extends results of Xiao--Zhu, Zhu--Zhuang, Zhou and Xie--Lu--Wang for triangulated, extriangulated, $(n+2)$-angulated and $n$-abelian categories,  respectively.\\[0.2cm]
\\[0.1cm]
\end{spacing}
\medskip
\end{abstract}

\pagestyle{myheadings}
\markboth{\rightline {\scriptsize J. He, J. Hu, D. Zhang and P. Zhou }}
         {\leftline{\scriptsize On the existence of Auslander-Reiten $n$-exangles in $n$-exangulated categories}}

\section{Introduction}
The notion of extriangulated categories was introduced by Nakaoka--Palu in \cite{NP}, which can be viewed
 as a simultaneous generalization of exact categories and triangulated categories.
 The data of such a category is a triplet $(\C,\E,\s)$, where $\C$ is an additive category, $\mathbb{E}: \C^{\rm op}\times \C \rightarrow \Ab$ is an additive bifunctor and $\mathfrak{s}$ assigns to each $\delta\in \mathbb{E}(C,A)$ a class of $3$-term sequences with end terms $A$ and $C$ such that certain axioms hold. Recently, Herschend--Liu--Nakaoka \cite{HLN}
introduced the notion of $n$-exangulated categories for any positive integer $n$. It is not only a higher dimensional analogue of extriangulated categories,
but also gives a common generalization of $n$-exact categories in the sense of
Jasso \cite{Ja} and $(n+2)$-angulated in the sense of Geiss--Keller--Oppermann \cite{GKO}.
 However,
 there are some other examples of $n$-exangulated categories which are neither $n$-exact nor $(n+2)$-angulated, see \cite{HLN, LZ,HZZ2}.

Auslander-Reiten theory was introduced by Auslander and Reiten in \cite{AR1,AR2}. Since its introduction, Auslander-Reiten theory has become a fundamental tool for
studying the representation theory of Artin algebras.
Later it has been generalized to these situation of exact categories \cite{Ji}, triangulated categories \cite{H,RV} and its subcategories \cite{AS,J}
and some certain additive categories \cite{L,J,S} by many authors.
Iyama, Nakaoka and Palu \cite{INP} developed  Auslander--Reiten theory for extriangulated categories.
This unifies Auslander--Reiten theories in exact categories and triangulated categories independently.
Xiao and Zhu \cite{XZ1,XZ2} showed that if a triangulated category $\C$ is locally finite, then $\C$ has Auslander-Reiten triangles. Recently, Zhu and Zhuang \cite{ZZ4} proved that if an extriangulated category $\C$ is locally finite, then $\C$ has Auslander-Reiten $\mathbb{E}$-triangles.
Later, Zhou \cite {Z2} extended Xiao-Zhu's result into $(n+2)$-angulated categories. Namely, Zhou proved that if an $(n+2)$-angulated  category $\C$ is locally finite, then $\C$ has Auslander-Reiten $(n+2)$-angules. Subsequently,  Xie-Lu-Wang \cite{XLW} proved a similar result to Zhou. More precisely, they showed that if an $n$-abelian category $\C$ is locally finite, then $\C$ has $n$-Auslander-Reiten sequences. Based on this idea, we have a natural question of whether their results of Zhou \cite{Z2} and Xie-Lu-Wang \cite{XLW} can be unified under the framework of $n$-exangulated categories or whether the result of Zhu-Zhuang \cite{ZZ4} has a higher counterpart. In this article, we give an affirmative answer.

Our main result is the following.

\begin{theorem}{\rm (see Theorem \ref{main} for details)}
Let $\C$ be a locally finite $n$-exangulated category. If $X\in\ind(\C)$ is a non-projective object, then there exists an Auslander-Reiten $n$-exangle ending at $X$, and if $Y\in\ind(\C)$ is a non-injective object, then there exists an Auslander-Reiten $n$-exangle starting at $Y$.
In this case, we say that $\C$ has Auslander-Reiten $n$-exangles.
\end{theorem}

This article is organized as follows: In Section 2, we recall the definition of $n$-exangulated category and review some results. In Section 3, we show our main result.

\section{Preliminaries}
In this section, we briefly review basic concepts and results concerning $n$-exangulated categories.

{ For any pair of objects $A,C\in\C$, an element $\del\in\E(C,A)$ is called an {\it $\E$-extension} or simply an {\it extension}. We also write such $\del$ as ${}_A\del_C$ when we indicate $A$ and $C$. The zero element ${}_A0_C=0\in\E(C,A)$ is called the {\it split $\E$-extension}. For any pair of $\E$-extensions ${}_A\del_C$ and ${}_{A'}\del{'}_{C'}$, let $\delta\oplus \delta'\in\mathbb{E}(C\oplus C', A\oplus A')$ be the
element corresponding to $(\delta,0,0,{\delta}{'})$ through the natural isomorphism $\mathbb{E}(C\oplus C', A\oplus A')\simeq\mathbb{E}(C, A)\oplus\mathbb{E}(C, A')
\oplus\mathbb{E}(C', A)\oplus\mathbb{E}(C', A')$.

For any $a\in\C(A,A')$ and $c\in\C(C',C)$,  $\E(C,a)(\del)\in\E(C,A')\ \ \text{and}\ \ \E(c,A)(\del)\in\E(C',A)$ are simply denoted by $a_{\ast}\del$ and $c^{\ast}\del$, respectively.

Let ${}_A\del_C$ and ${}_{A'}\del{'}_{C'}$ be any pair of $\E$-extensions. A {\it morphism} $(a,c)\colon\del\to{\delta}{'}$ of extensions is a pair of morphisms $a\in\C(A,A')$ and $c\in\C(C,C')$ in $\C$, satisfying the equality
$a_{\ast}\del=c^{\ast}{\delta}{'}$.}

\begin{definition}\cite[Definition 2.7]{HLN}
Let $\bf{C}_{\C}$ be the category of complexes in $\C$. As its full subcategory, define $\CC$ to be the category of complexes in $\C$ whose components are zero in the degrees outside of $\{0,1,\ldots,n+1\}$. Namely, an object in $\CC$ is a complex $X_{\bullet}=\{X_i,d^X_i\}$ of the form
\[ X_0\xrightarrow{d^X_0}X_1\xrightarrow{d^X_1}\cdots\xrightarrow{d^X_{n-1}}X_n\xrightarrow{d^X_n}X_{n+1}. \]
We write a morphism $f_{\bullet}\co X_{\bullet}\to Y_{\bullet}$ simply $f_{\bullet}=(f_0,f_1,\ldots,f_{n+1})$, only indicating the terms of degrees $0,\ldots,n+1$.
\end{definition}

\begin{definition}\cite[Definition 2.11]{HLN}
By Yoneda lemma, any extension $\del\in\E(C,A)$ induces natural transformations
\[ \del\ssh\colon\C(-,C)\Rightarrow\E(-,A)\ \ \text{and}\ \ \del\ush\colon\C(A,-)\Rightarrow\E(C,-). \]
For any $X\in\C$, these $(\del\ssh)_X$ and $\del\ush_X$ are given as follows.
\begin{enumerate}
\item[\rm(1)] $(\del\ssh)_X\colon\C(X,C)\to\E(X,A)\ :\ f\mapsto f\uas\del$.
\item[\rm (2)] $\del\ush_X\colon\C(A,X)\to\E(C,X)\ :\ g\mapsto g\sas\delta$.
\end{enumerate}
We simply denote $(\del\ssh)_X(f)$ and $\del\ush_X(g)$ by $\del\ssh(f)$ and $\del\ush(g)$, respectively.
\end{definition}

\begin{definition}\cite[Definition 2.9]{HLN}
 Let $\C,\E,n$ be as before. Define a category $\AE:=\AE^{n+2}_{(\C,\E)}$ as follows.
\begin{enumerate}
\item[\rm(1)]  A pair $\Xd$ is an object of the category $\AE$ with $X_{\bullet}\in\CC$
and $\del\in\E(X_{n+1},X_0)$, called an $\E$-attached
complex of length $n+2$, if it satisfies
$$(d_0^X)_{\ast}\del=0~~\textrm{and}~~(d^X_n)^{\ast}\del=0.$$
We also denote it by
$$X_0\xrightarrow{d_0^X}X_1\xrightarrow{d_1^X}\cdots\xrightarrow{d_{n-2}^X}X_{n-1}
\xrightarrow{d_{n-1}^X}X_n\xrightarrow{d_n^X}X_{n+1}\overset{\delta}{\dashrightarrow}.$$
\item[\rm (2)]  For such pairs $\Xd$ and $\langle Y_{\bullet},\rho\rangle$,  $f_{\bullet}\colon\Xd\to\langle Y_{\bullet},\rho\rangle$ is
defined to be a morphism in $\AE$ if it satisfies $(f_0)_{\ast}\del=(f_{n+1})^{\ast}\rho$.

\end{enumerate}
\end{definition}

\begin{definition}\cite[Definition 2.13]{HLN}\label{def1}
 An {\it $n$-exangle} is an object $\Xd$ in $\AE$ that satisfies the listed conditions.
\begin{enumerate}
\item[\rm (1)] The following sequence of functors $\C\op\to\Ab$ is exact.
$$
\C(-,X_0)\xrightarrow{\C(-,\ d^X_0)}\cdots\xrightarrow{\C(-,\ d^X_n)}\C(-,X_{n+1})\xrightarrow{~\del\ssh~}\E(-,X_0)
$$
\item[\rm (2)] The following sequence of functors $\C\to\Ab$ is exact.
$$
\C(X_{n+1},-)\xrightarrow{\C(d^X_n,\ -)}\cdots\xrightarrow{\C(d^X_0,\ -)}\C(X_0,-)\xrightarrow{~\del\ush~}\E(X_{n+1},-)
$$
\end{enumerate}
In particular any $n$-exangle is an object in $\AE$.
A {\it morphism of $n$-exangles} simply means a morphism in $\AE$. Thus $n$-exangles form a full subcategory of $\AE$.
\end{definition}

\begin{definition}\cite[Definition 2.22]{HLN}
Let $\s$ be a correspondence which associates a homotopic equivalence class $\s(\del)=[{}_A{X_{\bullet}}_C]$ to each extension $\del={}_A\del_C$. Such $\s$ is called a {\it realization} of $\E$ if it satisfies the following condition for any $\s(\del)=[X_{\bullet}]$ and any $\s(\rho)=[Y_{\bullet}]$.
\begin{itemize}
\item[{\rm (R0)}] For any morphism of extensions $(a,c)\co\del\to\rho$, there exists a morphism $f_{\bullet}\in\CC(X_{\bullet},Y_{\bullet})$ of the form $f_{\bullet}=(a,f_1,\ldots,f_n,c)$. Such $f_{\bullet}$ is called a {\it lift} of $(a,c)$.
\end{itemize}
In such a case, we simple say that \lq\lq$X_{\bullet}$ realizes $\del$" whenever they satisfy $\s(\del)=[X_{\bullet}]$.

Moreover, a realization $\s$ of $\E$ is said to be {\it exact} if it satisfies the following conditions.
\begin{itemize}
\item[{\rm (R1)}] For any $\s(\del)=[X_{\bullet}]$, the pair $\Xd$ is an $n$-exangle.
\item[{\rm (R2)}] For any $A\in\C$, the zero element ${}_A0_0=0\in\E(0,A)$ satisfies
\[ \s({}_A0_0)=[A\ov{\id_A}{\lra}A\to0\to\cdots\to0\to0]. \]
Dually, $\s({}_00_A)=[0\to0\to\cdots\to0\to A\ov{\id_A}{\lra}A]$ holds for any $A\in\C$.
\end{itemize}
Note that the above condition {\rm (R1)} does not depend on representatives of the class $[X_{\bullet}]$.
\end{definition}

\begin{definition}\cite[Definition 2.23]{HLN}
Let $\s$ be an exact realization of $\E$.
\begin{enumerate}
\item[\rm (1)] An $n$-exangle $\Xd$ is called an $\s$-{\it distinguished} $n$-exangle if it satisfies $\s(\del)=[X_{\bullet}]$. We often simply say {\it distinguished $n$-exangle} when $\s$ is clear from the context.
\item[\rm (2)]  An object $X_{\bullet}\in\CC$ is called an {\it $\s$-conflation} or simply a {\it conflation} if it realizes some extension $\del\in\E(X_{n+1},X_0)$.
\item[\rm (3)]  A morphism $f$ in $\C$ is called an {\it $\s$-inflation} or simply an {\it inflation} if it admits some conflation $X_{\bullet}\in\CC$ satisfying $d_0^X=f$.
\item[\rm (4)]  A morphism $g$ in $\C$ is called an {\it $\s$-deflation} or simply a {\it deflation} if it admits some conflation $X_{\bullet}\in\CC$ satisfying $d_n^X=g$.
\end{enumerate}
\end{definition}

\begin{definition}\cite[Definition 2.27]{HLN}
For a morphism $f_{\bullet}\in\CC(X_{\bullet},Y_{\bullet})$ satisfying $f_0=\id_A$ for some $A=X_0=Y_0$, its {\it mapping cone} $M_{_{\bullet}}^f\in\CC$ is defined to be the complex
\[ X_1\xrightarrow{d^{M_f}_0}X_2\oplus Y_1\xrightarrow{d^{M_f}_1}X_3\oplus Y_2\xrightarrow{d^{M_f}_2}\cdots\xrightarrow{d^{M_f}_{n-1}}X_{n+1}\oplus Y_n\xrightarrow{d^{M_f}_n}Y_{n+1} \]
where $d^{M_f}_0=\begin{bmatrix}-d^X_1\\ f_1\end{bmatrix},$
$d^{M_f}_i=\begin{bmatrix}-d^X_{i+1}&0\\ f_{i+1}&d^Y_i\end{bmatrix}\ (1\le i\le n-1),$
$d^{M_f}_n=\begin{bmatrix}f_{n+1}&d^Y_n\end{bmatrix}$.

{\it The mapping cocone} is defined dually, for morphisms $h_{\bullet}$ in $\CC$ satisfying $h_{n+1}=\id$.
\end{definition}

\begin{definition}\cite[Definition 2.32]{HLN}
An {\it $n$-exangulated category} is a triplet $(\C,\E,\s)$ of additive category $\C$, additive bifunctor $\E\co\C\op\times\C\to\Ab$, and its exact realization $\s$, satisfying the following conditions.
\begin{itemize}
\item[{\rm (EA1)}] Let $A\ov{f}{\lra}B\ov{g}{\lra}C$ be any sequence of morphisms in $\C$. If both $f$ and $g$ are inflations, then so is $g\circ f$. Dually, if $f$ and $g$ are deflations, then so is $g\circ f$.

\item[{\rm (EA2)}] For $\rho\in\E(D,A)$ and $c\in\C(C,D)$, let ${}_A\langle X_{\bullet},c\uas\rho\rangle_C$ and ${}_A\langle Y_{\bullet},\rho\rangle_D$ be distinguished $n$-exangles. Then $(\id_A,c)$ has a {\it good lift} $f_{\bullet}$, in the sense that its mapping cone gives a distinguished $n$-exangle $\langle M^f_{_{\bullet}},(d^X_0)\sas\rho\rangle$.
 \item[{\rm (EA2$\op$)}] Dual of {\rm (EA2)}.
\end{itemize}
Note that the case $n=1$, a triplet $(\C,\E,\s)$ is a  $1$-exangulated category if and only if it is an extriangulated category, see \cite[Proposition 4.3]{HLN}.
\end{definition}

\begin{example}
From \cite[Proposition 4.34]{HLN} and \cite[Proposition 4.5]{HLN},  we know that $n$-exact categories and $(n+2)$-angulated categories are $n$-exangulated categories.
There are some other examples of $n$-exangulated categories
 which are neither $n$-exact nor $(n+2)$-angulated, see \cite{HLN,LZ,HZZ2}.
\end{example}
The following some Lemmas are very useful which are needed later on.
\begin{lemma}\emph{\cite[Proposition 3.6]{HLN}}\label{a2}
\rm Let ${}_A\langle X_{\bullet},\delta\rangle_C$ and ${}_B\langle Y_{\bullet},\rho\rangle_D$ be distinguished $n$-exangles. Suppose that we are given a commutative square
$$\xymatrix{
 X_0 \ar[r]^{{d_0^X}} \ar@{}[dr]|{\circlearrowright} \ar[d]_{a} & X_1 \ar[d]^{b}\\
 Y_0  \ar[r]_{d_0^Y} &Y_1
}
$$
in $\C$. Then there is a morphism $f_{\bullet}\colon \langle X_{\bullet},\delta\rangle\to\langle Y_{\bullet},\rho\rangle$ which satisfies $f_0=a$ and $f_1=b$.
\end{lemma}
\begin{lemma}\label{y1}
Let $\C$ be an $n$-exangulated category , and $$\xymatrix{
X_0\ar[r]^{f_0}\ar@{}[dr] \ar[d]^{a_0} &X_1 \ar[r]^{f_1} \ar@{}[dr]\ar[d]^{a_1}\ar@{-->}[dl]^{h_1} &X_2 \ar[r]^{f_2} \ar@{}[dr]\ar[d]^{a_2}\ar@{-->}[dl]^{h_2}&\cdot\cdot\cdot \ar[r]\ar@{}[dr] &X_n \ar[r]^{f_n} \ar@{}[dr]\ar[d]^{a_n}&X_{n+1} \ar@{}[dr]\ar[d]^{a_{n+1}} \ar@{-->}[dl]^{h_{n+1}}\ar@{-->}[r]^-{\delta} &\\
{Y_0}\ar[r]^{g_0} &{Y_1}\ar[r]^{g_1}&{Y_2} \ar[r]^{g_2} &\cdot\cdot\cdot \ar[r] &{Y _n}\ar[r]^{g_n}  &{Y_{n+1}} \ar@{-->}[r]^-{\eta} &}
$$
any morphism of distinguished $n$-exangles. Then the following are equivalent:
\begin{itemize}
\item[\rm (1)]There is a morphism $h_1\colon X_1\to Y_0$, such that $h_1f_0=a_0$.

\item[\rm (2)]There is a morphism $h_{n+1}\colon X_{n+1}\to Y_n$, such that $g_nh_{n+1}=a_{n+1}$.

\item[\rm (3)] $ (a_0)_{*}{\delta}=(a_{n+1})^{*}{\eta}=0$.

\item[\rm (4)] $a_{\bullet}=(a_0,a_1,\cdot\cdot\cdot,a_{n+1})\colon\Xd\to\langle Y_{\bullet},\eta\rangle$ is null-homotopic.

\end{itemize}
\end{lemma}
\proof This follows from \cite[Lemma 3.3]{ZW} and \cite[Lemma 3.6]{F}.
\qed

\begin{corollary}\label{y2}
If $a^{\mr}$ is the identity on $\Xd$ as above,  then the following are equivalent:
\begin{itemize}
\item[\rm (1)] $f_0$ is a split monomorphism (also known as a section).

\item[\rm (2)] $f_n$ is a split epimorphism (also known as a retraction).

\item[\rm (3)] $ {\delta}=0$.

\item[\rm (4)] $a_{\bullet}$ is null-homotopic.

\end{itemize}
If a distinguished $n$-exangle satisfies one of the above equivalent conditions, it is called \emph{split}.
\end{corollary}
\begin{definition}\label{def2}\cite[Definition 3.14 ]{ZW}and \cite[Definition 3.2]{LZ}
Let $(\C,\E,\s)$ be an $n$-exangulated category. An object $P\in\C$ is called \emph{projective} if, for any distinguished $n$-exangle
$$A_0\xrightarrow{\alpha_0}A_1\xrightarrow{\alpha_1}A_2\xrightarrow{\alpha_2}\cdots\xrightarrow{\alpha_{n-2}}A_{n-1}
\xrightarrow{\alpha_{n-1}}A_n\xrightarrow{\alpha_n}A_{n+1}\overset{\delta}{\dashrightarrow}$$
and any morphism $c$ in $\C(P,A_{n+1})$, there exists a morphism $b\in\C(P,A_n)$ satisfying $\alpha_n\circ b=c$.
We denote the full subcategory of projective objects in $\C$ by $\P$.
Dually, the full subcategory of injective objects in $\C$ is denoted by $\I$.
\end{definition}

\begin{lemma}\rm\label{def2}\cite[Lemma 3.4 ]{LZ} Let $(\C,\E,\s)$ be an $n$-exangulated category. Then the following statements are equivalent for an object $P\in\C$.
\begin{itemize}
\item[\rm (1)] $\E(P,A)=0$ for any $A\in\C$.
\item[\rm (2)] $P$ is projective.
\item[\rm (3)] Any distinguished $n$-exangle
$A_0\xrightarrow{\alpha_0}A_1\xrightarrow{\alpha_1}A_2\xrightarrow{\alpha_2}\cdots\xrightarrow{\alpha_{n-2}}A_{n-1}
\xrightarrow{\alpha_{n-1}}A_n\xrightarrow{\alpha_n}P\overset{\delta}{\dashrightarrow}$ splits.
\end{itemize}
\end{lemma}

We denote by ${\rm rad}_{\C}$ the Jacobson radical of $\C$. Namely, ${\rm rad}_{\C}$ is an ideal of $\C$ such that ${\rm rad}_{\C}(A, A)$
coincides with the Jacobson radical of the endomorphism ring ${\rm End}(A)$ for any $A\in\C$.

\begin{definition}\cite[Definition 3.3 ]{HZ} When $n\geq2$, a distinguished $n$-exangle in $\C$ of the form
$$A_{\bullet}:~~A_0\xrightarrow{\alpha_0}A_1\xrightarrow{\alpha_1}A_2\xrightarrow{\alpha_2}\cdots\xrightarrow{\alpha_{n-2}}A_{n-1}
\xrightarrow{\alpha_{n-1}}A_n\xrightarrow{\alpha_n}A_{n+1}\overset{}{\dashrightarrow}$$
is minimal if $\alpha_1,\alpha_2,\cdots,\alpha_{n-1}$ are in $\rad_{\C}$.

\end{definition}
The following lemma shows that $\E$-extension in an equivalence class can be chosen in a minimal way in a Krull-Schmidt $n$-exangulated category.

\begin{lemma}\rm\label{ml}\cite[Lemma 3.4 ]{HZ}  Let $\C$ be a Krull-Schmidt $n$-exangulated category, $A_0,A_{n+1}\in\C$. Then for every equivalence class of $\E$-extension of $A_0$ by $A_{n+1}$, there is a representation
$$A_{\bullet}:~~A_0\xrightarrow{\alpha_0}A_1\xrightarrow{\alpha_1}A_2\xrightarrow{\alpha_2}\cdots\xrightarrow{\alpha_{n-2}}A_{n-1}
\xrightarrow{\alpha_{n-1}}A_n\xrightarrow{\alpha_n}A_{n+1}\overset{}{\dashrightarrow}$$
such that $\alpha_1,\alpha_2,\cdots,\alpha_{n-1}$ are in $\rad_{\C}$. Moreover, $A_{\bullet}$ is a direct summand of every other equivalent $\E$-extension.
\end{lemma}

\begin{remark} Let $\C$ be a Krull-Schmidt $n$-exangulated category. By the Krull-Schmidt property of $\C$, every minimal distinguished $n$-exangle in each equivalence class is unique up to isomorphism.
\end{remark}
\section{Locally finite $n$-exangulated categories }
In this section, let $k$ be a field. We always assume that $\C$ is a $k$-linear Hom-finite Krull-Schmidt $n$-exangulated category.
We denote by $\ind(\C)$ the set of isomorphism classes of indecomposable objects in $\C$.
\begin{definition}\label{222} Let $\C$ be an $n$-exangulated category.
A distinguished $n$-exangle
$$A_0\xrightarrow{\alpha_0}A_1\xrightarrow{\alpha_1}A_2\xrightarrow{\alpha_2}\cdots\xrightarrow{\alpha_{n-2}}A_{n-1}
\xrightarrow{\alpha_{n-1}}A_n\xrightarrow{\alpha_n}A_{n+1}\overset{\delta}{\dashrightarrow}$$
in $\C$ is called an \emph{Auslander-Reiten $n$-exangle }if
$\alpha_0$ is left almost split, $\alpha_n$ is right almost split and
when $n\geq 2$, $\alpha_1,\alpha_2,\cdots,\alpha_{n-1}$ are in $\rad_{\C}$.
\end{definition}
\begin{remark} (1) If $\C$ is an $n$-abelian category, then Definition \ref{222} coincides with the definition of $n$-Auslander-Reiten sequence of $n$-abelian category (cf. \cite{XLW}).

(2) If $\C$ is an $({n+2})$-angulated category, then Definition \ref{222} coincides with the definition of Auslander-Reiten $({n+2})$-angle of $({n+2})$-angulated category (cf. \cite{Z2}).
\end{remark}

\begin{lemma}\label{r1}
Let $\C$ be an $n$-exangulated category and
$$A_{\bullet}:~~A_0\xrightarrow{\alpha_0}A_1\xrightarrow{\alpha_1}A_2\xrightarrow{\alpha_2}\cdots\xrightarrow{\alpha_{n-2}}A_{n-1}
\xrightarrow{\alpha_{n-1}}A_n\xrightarrow{\alpha_n}A_{n+1}\overset{\delta}{\dashrightarrow}$$
be a distinguished $n$-exangle in $\C$. Then the following statements are equivalent:
\begin{itemize}
\item[\rm (1)] $A_{\bullet}$ is an Auslander-Reiten $n$-exangle;
\item[\rm (2)] ${\rm{End}}(A_0)$ is local, if $n\geq 2$, $ \alpha_1,\cdots,\alpha_{n-1}$ are in ${\rm rad}_{\C}$ and $\alpha_n$ is right almost split;
\item[\rm (3)] ${\rm{End}}(A_{n+1})$ is local, if $n\geq 2$, $\alpha_1,\alpha_2,\cdots,\alpha_{n-1}$ are in ${\rm rad}_{\C}$ and $\alpha_0$ is left almost split.
\end{itemize}
\end{lemma}

\proof The proof given in \cite[Lemma 5.3]{F1} can be adapted to the context of $n$-exangulated
categories, we omit it.  \qed
\medskip

For any $X\in\ind(\C)$, we denote by $\supp \Hom_{\C}(X,-)$ the subcategory of $\C$ generated by
objects $Y$ in $\ind(\C)$ with $\Hom_{\C}(X,Y)\neq0$. Similarly, $\supp\Hom_{\C}(-,X)$ denotes
the subcategory generated by objects $Y$ in $\ind(\C)$ with $\Hom_{\C}(Y,X)\neq 0$. If
$\supp\Hom_{\C}(X,-)$ ($\supp\Hom_{\C}(-,X)$, respectively) contains only finitely many
indecomposables, we say that $|\supp\Hom_{\C}(X,-)|<\infty $ ($|\supp\Hom_{\C}(-,X)|<\infty$
respectively).
\medskip

Based on the definition of locally finite $(n+2)$-angulated categories and locally finite $n$-abelian categories, \cite{Z2,XLW},
we define the notion of locally finite $n$-exangulated categories.
\begin{definition}
A  $n$-exangulated category $\C$ is called \emph{locally finite}
if $|\supp\Hom_{\C}(X,-)|<\infty$ and $|\supp\Hom_{\C}(-,X)|<\infty$, for any object $X\in\ind(\C)$.
\end{definition}
\begin{definition}\label{y2}
Let $\C$ be an $n$-exangulated category and $X_{n+1}, Y_{0} \in\ind(\C)$. We define a set of distinguished $n$-exangles as follows:
{\small $$
S(X_{n+1}):=\biggl\{X_{\bullet}:X_0\xrightarrow{\alpha_0}\cdots\xrightarrow{\alpha_{n-1}}X_n\xrightarrow{\alpha_n}X_{n+1}\overset{\delta}{\dashrightarrow}~\bigg|~\begin{array}{ll}
 X_{\bullet}~\mbox{is a non-split distinguished $n$-exangle}\\
\textrm{with}~X_0\in\ind(\C), \mbox{and when}\\ n\geq2,~\alpha_1,\alpha_2,\cdots,\alpha_{n-1}$ in ${\rm rad_{\C}}.
 \end{array}\biggl\}$$}

Dually, we can define a set of distinguished $n$-exangles as follows:
$$
T(Y_{0}):=\biggl\{Y_{\bullet}:Y_0\xrightarrow{\beta_0}\cdots\xrightarrow{\beta_{n-1}}Y_n\xrightarrow{\beta_n}Y_{n+1}\overset{\eta}{\dashrightarrow}~\bigg|~\begin{array}{ll}
 Y_{\bullet}~\mbox{is a non-split distinguished $n$-exangle}\\
\textrm{with}~Y_{n+1}\in\ind(\C), \mbox{and when}\\ n\geq2,~ \beta_1,\beta_2,\cdots,\beta_{n-1}$ in ${\rm rad_{\C}}.
 \end{array}\biggl\}$$
\end{definition}

\begin{lemma}\rm\label{lem1}Let $(\C,\E,\s)$ be an $n$-exangulated category.
\begin{itemize}
\item[ (1)] If $X_{n+1}\in\ind(\C)$ is a non-projective object, then $S(X_{n+1})$ is non-empty.
\item[ (2)] If $Y_{0}\in\ind(\C)$ is a non-injective object, then $T(Y_{0})$ is non-empty.

\end{itemize}
\end{lemma}
\proof We only show that $(1)$, dually one can prove $(2)$.

Since $X_{n+1}\in\ind(\C)$ is a non-projective, there is an object $X_{0}\in\C$, such that $\E(X_{n+1},X_{0})\neq 0$ by Lemma \ref{def2}. That is to say, there exists a non-split distinguished $n$-exangle:$$X_{\bullet}:X_0\xrightarrow{\alpha_0}X_1\xrightarrow{\alpha_1}X_2\xrightarrow{\alpha_2}\cdots\xrightarrow{\alpha_{n-1}}X_n\xrightarrow{\alpha_n}X_{n+1}\overset{\delta}{\dashrightarrow}.$$
Since $\C$ is a Krull-Schmidt category, we decompose $X_0$ into a direct sum of indecomposable objects $X_0=\bigoplus\limits_{i=1}^dA_i$.
Without loss of generality, we can assume that $X_0=U\oplus V$ where $U$ and $V$ are indecomposable. Since $\E(X_{n+1},X_{0})\simeq\mathbb{E}(X_{n+1}, U\oplus V)\simeq\mathbb{E}(X_{n+1}, U)\oplus\mathbb{E}(X_{n+1}, V)$.
 We claim that the at least one of the following two distinguished $n$-exangles is non-split
$$U\xrightarrow{\gamma_0}C_1\xrightarrow{\gamma_1}C_2\xrightarrow{\gamma_2}\cdots\xrightarrow{\gamma_{n-1}}C_n\xrightarrow{\gamma_n}X_{n+1}\overset{\eta}{\dashrightarrow}$$
$$V\xrightarrow{\beta_0}D_1\xrightarrow{\beta_1}D_2\xrightarrow{\beta_2}\cdots\xrightarrow{\beta_{n-1}}D_n\xrightarrow{\beta_n}X_{n+1}\overset{\eta^{'}}{\dashrightarrow}$$
where $\eta:=\left(\begin{smallmatrix}1&0\end{smallmatrix}\right)_*\delta$ and $\eta':=\left(\begin{smallmatrix}0&1\end{smallmatrix}\right)_*\delta$. Otherwise, $\delta=\eta\oplus\eta'=0\in\E(X_{n+1},X_{0})$. This is a contradiction  since $\delta$ is non-split.

We can take an distinguished $n$-exangle as we want by Lemma \ref{ml}. This completes the proof.\qed

\begin{definition}\label{y6}
Let $\C$ be an $n$-exangulated category, and
$$X_{\bullet}:X_0\xrightarrow{\alpha_0}X_1\xrightarrow{\alpha_1}X_2\xrightarrow{\alpha_2}\cdots\xrightarrow{\alpha_{n-1}}X_n\xrightarrow{\alpha_n}X_{n+1}\overset{\delta}{\dashrightarrow}$$
$$U_{\bullet}:U_0\xrightarrow{\beta_0}U_1\xrightarrow{\beta_1}U_2\xrightarrow{\beta_2}\cdots\xrightarrow{\beta_{n-1}}U_n\xrightarrow{\beta_n}X_{n+1}\overset{\delta^{'}}{\dashrightarrow}$$
two distinguished $n$-exangles in $S(X_{n+1})$. We say that
$X_{\bullet}>U_{\bullet}$ if there exists a morphism of distinguished $n$-exangles as follows:
$$\xymatrix{
X_0\ar[r]^{\alpha_0}\ar@{}[dr] \ar@{-->}[d]^{\varphi_0} &X_1 \ar[r]^{\alpha_1} \ar@{}[dr]\ar@{-->}[d]^{\varphi_1}&X_2 \ar[r]^{\alpha_2} \ar@{}[dr]\ar@{-->}[d]^{\varphi_2}&\cdot\cdot\cdot \ar[r]\ar@{}[dr] &X_n \ar[r]^{\alpha_n} \ar@{}[dr]\ar@{-->}[d]^{\varphi_n}&X_{n+1} \ar@{}[dr]\ar@{=}[d]^{} \ar@{-->}[r]^-{\delta} &\\
{U_0}\ar[r]^{\beta_0} &{U_1}\ar[r]^{\beta_1}&{U_2} \ar[r]^{\beta_2} &\cdot\cdot\cdot \ar[r] &{U _n}\ar[r]^{\beta_n}  &{X_{n+1}} \ar@{-->}[r]^-{\delta^{'}} &}
$$
We say that $X_{\bullet}\sim U_{\bullet}$ if $\varphi_0$ is an isomorphism.

 Dually, let
$$Y_{\bullet}:Y_0\xrightarrow{\alpha_0}Y_1\xrightarrow{\alpha_1}Y_2\xrightarrow{\alpha_2}\cdots\xrightarrow{\alpha_{n-1}}Y_n\xrightarrow{\alpha_n}Y_{n+1}\overset{\eta}{\dashrightarrow}$$
$$V_{\bullet}:Y_0\xrightarrow{\beta_0}V_1\xrightarrow{\beta_1}V_2\xrightarrow{\beta_2}\cdots\xrightarrow{\beta_{n-1}}V_n\xrightarrow{\beta_n}V_{n+1}\overset{\eta^{'}}{\dashrightarrow}$$
two distinguished $n$-exangles in $T(Y_{0})$. We say that
$Y_{\bullet}>V_{\bullet}$ if there exists a morphism of distinguished $n$-exangles as follows:
$$\xymatrix{
Y_0\ar[r]^{\alpha_0}\ar@{}[dr] \ar@{=}[d]^{\varphi_0} &Y_1 \ar[r]^{\alpha_1} \ar@{}[dr]\ar@{-->}[d]^{\varphi_1}&Y_2 \ar[r]^{\alpha_2} \ar@{}[dr]\ar@{-->}[d]^{\varphi_2}&\cdot\cdot\cdot \ar[r]\ar@{}[dr] &Y_n \ar[r]^{\alpha_n} \ar@{}[dr]\ar@{-->}[d]^{\varphi_n}&Y_{n+1} \ar@{}[dr]\ar@{-->}[d]^{\varphi_{n+1}} \ar@{-->}[r]^-{\eta} &\\
{Y_0}\ar[r]^{\beta_0} &{V_1}\ar[r]^{\beta_1}&{V_2} \ar[r]^{\beta_2} &\cdot\cdot\cdot \ar[r] &{V_n}\ar[r]^{\beta_n}  &{V_{n+1}} \ar@{-->}[r]^-{\eta^{'}} &}
$$
We say that $Y_{\bullet}\sim V_{\bullet}$ if $\varphi_{n+1}$ is an isomorphism.

\end{definition}

In the following, we will consider a direct ordered set, namely, a partially ordered set with every pair of elements has a lower bound.

\begin{lemma}\label{y7}
$S(X_{n+1})$ is a direct ordered set with the relation defined in Definition \ref{y6}, and
$T(Y_{0})$ is a direct ordered set with the relation defined in Definition \ref{y6}.
\end{lemma}
\proof We just prove the first statement, the second statement proves similarly.

Assume that
$$X_{\bullet}:X_0\xrightarrow{\alpha_0}X_1\xrightarrow{\alpha_1}X_2\xrightarrow{\alpha_2}\cdots\xrightarrow{\alpha_{n-1}}X_n\xrightarrow{\alpha_n}X_{n+1}\overset{\delta}{\dashrightarrow}$$
and $$U_{\bullet}:U_0\xrightarrow{\beta_0}U_1\xrightarrow{\beta_1}U_2\xrightarrow{\beta_2}\cdots\xrightarrow{\beta_{n-1}}U_n\xrightarrow{\beta_n}X_{n+1}\overset{\delta^{'}}{\dashrightarrow}$$
belong to $S(X_{n+1})$.

Firstly, the axioms of reflexivity and transitivity are clear. Secondly, We show that if $X_{\bullet}>U_{\bullet}$ and $U_{\bullet}>X_{\bullet}$, then
$X_{\bullet}\sim U_{\bullet}$.

Since $X_{\bullet}>U_{\bullet}$ and $U_{\bullet}>X_{\bullet}$, we have the following two commutative diagrams
$$\xymatrix{
X_0\ar[r]^{\alpha_0}\ar@{}[dr] \ar@{-->}[d]^{\varphi_0} &X_1 \ar[r]^{\alpha_1} \ar@{}[dr]\ar@{-->}[d]^{\varphi_1}&X_2 \ar[r]^{\alpha_2} \ar@{}[dr]\ar@{-->}[d]^{\varphi_2}&\cdot\cdot\cdot \ar[r]\ar@{}[dr] &X_n \ar[r]^{\alpha_n} \ar@{}[dr]\ar@{-->}[d]^{\varphi_n}&X_{n+1} \ar@{}[dr]\ar@{=}[d]^{} \ar@{-->}[r]^-{\delta} &\\
{U_0}\ar[r]^{\beta_0} &{U_1}\ar[r]^{\beta_1}&{U_2} \ar[r]^{\beta_2} &\cdot\cdot\cdot \ar[r] &{U _n}\ar[r]^{\beta_n}  &{X_{n+1}} \ar@{-->}[r]^-{\delta^{'}} &}
$$
$$\xymatrix{
U_0\ar[r]^{\beta_0}\ar@{}[dr] \ar@{-->}[d]^{\psi_0} &U_1 \ar[r]^{\beta_1} \ar@{}[dr]\ar@{-->}[d]^{\psi_1}&U_2 \ar[r]^{\beta_2} \ar@{}[dr]\ar@{-->}[d]^{\psi_2}&\cdot\cdot\cdot \ar[r]\ar@{}[dr] &U_n \ar[r]^{\beta_n} \ar@{}[dr]\ar@{-->}[d]^{\psi_n}&X_{n+1} \ar@{}[dr]\ar@{=}[d]^{} \ar@{-->}[r]^-{\delta^{'}} &\\
{X_0}\ar[r]^{\alpha_0} &{X_1}\ar[r]^{\alpha_1}&{X_2} \ar[r]^{\alpha_2} &\cdot\cdot\cdot \ar[r] &{X _n}\ar[r]^{\alpha_n}  &{X_{n+1}} \ar@{-->}[r]^-{\delta} &}
$$
We claim that $\psi_0\varphi_0$ is an isomorphism. Since $X_0$ is an indecomposable, we have that $\End(X_0)$ is local implies that
$\psi_0\varphi_0$ is nilpotent or is
an isomorphism. If $\psi_0\varphi_0$ is nilpotent, there exists a positive integer $m$ such that $(\psi_0\varphi_0)^m=0$. We write $\omega_i=\psi_i\varphi_i, i=1,2,\cdot\cdot\cdot,n$.
Thus we have the following commutative diagram
$$\xymatrix{
X_0\ar[r]^{\alpha_0}\ar@{}[dr] \ar@{-->}[d]^{(\psi_0\varphi_0)^m} &X_1 \ar[r]^{\alpha_1} \ar@{}[dr]\ar@{-->}[d]^{\omega_1^m}&X_2 \ar[r]^{\alpha_2} \ar@{}[dr]\ar@{-->}[d]^{\omega_2^m}&\cdot\cdot\cdot \ar[r]\ar@{}[dr] &X_n \ar[r]^{\alpha_n} \ar@{}[dr]\ar@{-->}[d]^{\omega_n^m}&X_{n+1} \ar@{}[dr]\ar@{=}[d]^{} \ar@{-->}[r]^-{\delta} &\\
{X_0}\ar[r]^{\alpha_0} &{X_1}\ar[r]^{\alpha_1}&{X_2} \ar[r]^{\alpha_2} &\cdot\cdot\cdot \ar[r] &{X _n}\ar[r]^{\alpha_n}  &{X_{n+1}} \ar@{-->}[r]^-{\delta} &}
$$
Then $\delta={(\psi_0\varphi_0)_{\ast}^m}\delta=0$. This is a contradiction by Corollary \ref{y2} since $X_{\bullet}$ is non-split.
Hence $\psi_0\varphi_0$ is an isomorphism.
By a similar argument we obtain that $\varphi_{0}\psi_{0}$ is an isomorphism.
This shows that $\varphi_0$ is isomorphism. So $X_{\bullet}\sim U_{\bullet}$.

Finally, we show that if $X_{\bullet},U_{\bullet}\in S(X_{n+1})$, then there exists $C_{\bullet}\in S(X_{n+1})$ such that $X_{\bullet}>C_{\bullet}$ and
$U_{\bullet}> C_{\bullet}$.

For the morphism $\beta_{n}\colon U_n\to X_{n+1}$, by {\rm (EA2)}, we can observe that $(\id_{X_{0}},\beta_n)$ has a {\it good lift} $f^{\mr}=(\id_{X_{0}},\psi_1,\ldots,\psi_n,\beta_n)$, that is,
there exists the following commutative diagram of distinguished $n$-exangles
$$\xymatrix{
X_0\ar[r]^{\gamma_0}\ar@{}[dr] \ar@{=}[d]^{} &Z_1 \ar[r]^{\gamma_1} \ar@{}[dr]\ar@{-->}[d]^{\psi_1}&Z_2 \ar[r]^{\gamma_2} \ar@{}[dr]\ar@{-->}[d]^{\psi_2}&\cdot\cdot\cdot \ar[r]\ar@{}[dr] &Z_n \ar[r]^{\gamma_n} \ar@{}[dr]\ar@{-->}[d]^{\psi_n}&U_{n} \ar@{}[dr]\ar[d]^{\beta_{n}} \ar@{-->}[r]^-{\beta^{\ast}_{n}\delta} &\\
{X_0}\ar[r]^{\alpha_0} &{X_1}\ar[r]^{\alpha_1}&{X_2} \ar[r]^{\alpha_2} &\cdot\cdot\cdot \ar[r] &{X _n}\ar[r]^{\alpha_n}  &{X_{n+1}} \ar@{-->}[r]^-{\delta} &}
$$
such that
$M_{\bullet}:~Z_1\xrightarrow{} M_1\xrightarrow{}
 M_2\xrightarrow{}\cdots\xrightarrow{}M_{n-1}
                                     \xrightarrow{}U_n\oplus X_n\xrightarrow{\left(
                                       \begin{smallmatrix}
                                         \beta_n,&\alpha_{n}
                                       \end{smallmatrix}
                                     \right)}X_{n+1}\overset{(\gamma_0)_{\ast}\delta}{\dashrightarrow}$
 is a distinguished $n$-exangle in $\C$, where $M_i=Z_{i+1}\oplus X_{i}~(i=1,2,\cdots,n-1)$.
Since $\beta_n$ and $\alpha_n$ are not split epimorphisms, we have that
$(\beta_n,\alpha_n)$ is also not split epimorphism.
Otherwise, there exists a morphism $\binom{s}{t}\colon X_{n+1}\to U_n\oplus X_n$
such that $(\beta_n,\alpha_n)\binom{s}{t}=1_{X_{n+1}}$ and then $\beta_ns+\alpha_nt=1_{X_{n+1}}$.
Since $X_{n+1}$ is an indecomposable, we have that $\End(X_{n+1})$ is local implies that
either $\beta_ns$ or $\alpha_nt$ is an isomorphism.
Thus either $\beta_n$ or $\alpha_n$ is a split epimorphism, a contradiction.
That is, $M_{\bullet}$ is non-split.

Without loss of generality, we can assume that $Z_1=U\oplus V$ where $U$ and $V$ are indecomposable. For the morphism $p_1=(1,0)\colon U\oplus V\to U$, by  (EA2$\op$), we can observe that $(p_1,\id_{X_{n+1}})$ has a {\it good lift} $g^{\mr}=(p_1,\varphi_1,\ldots,\varphi_n,\id_{X_{n+1}})$, that is,
there exists the following commutative diagram of distinguished $n$-exangles
$$\xymatrix{
U\oplus V\ar[r]^{\;\;(u,~v)}\ar@{}[dr] \ar[d]^{p_1} &M_1 \ar[r]^{} \ar@{}[dr]\ar@{-->}[d]^{\varphi_1}&M_2 \ar[r]^{} \ar@{}[dr]\ar@{-->}[d]^{}&\cdot\cdot\cdot \ar[r]\ar@{}[dr] &M_{n-1} \ar[r]^{} \ar@{}[dr]\ar@{-->}[d]^{}&U_n\oplus X_n\ar[r]^{} \ar@{}[dr]\ar@{-->}[d]^{}&{X_{n+1}} \ar@{}[dr]\ar@{=}[d]^{} \ar@{-->}[r]^-{} &\\
{U}\ar[r]^{\delta_0} &{L_1}\ar[r]^{}&{L_2} \ar[r]^{} &\cdot\cdot\cdot \ar[r] &{L_{n-1}}\ar[r]^{}  &{L_{n}}\ar[r]^{}  &{X_{n+1}} \ar@{-->}[r]^-{} &}
$$
Similarly, for the morphism $p_2=(0,1)\colon U\oplus V\to V$, there exists
the following commutative diagram of distinguished $n$-exangles
\vspace{-2mm}
$$\xymatrix{
U\oplus V\ar[r]^{\;\;(u,~v)}\ar@{}[dr] \ar[d]^{p_2} &M_1 \ar[r]^{} \ar@{}[dr]\ar@{-->}[d]^{m_1}&M_2 \ar[r]^{} \ar@{}[dr]\ar@{-->}[d]^{}&\cdot\cdot\cdot \ar[r]\ar@{}[dr] &M_{n-1} \ar[r]^{} \ar@{}[dr]\ar@{-->}[d]^{}&U_n\oplus X_n\ar[r]^{} \ar@{}[dr]\ar@{-->}[d]^{}&{X_{n+1}} \ar@{}[dr]\ar@{=}[d]^{} \ar@{-->}[r]^-{} &\\
{V}\ar[r]^{\eta_0} &{N_1}\ar[r]^{}&{N_2} \ar[r]^{} &\cdot\cdot\cdot \ar[r] &{N_{n-1}}\ar[r]^{}  &{N_{n}}\ar[r]^{}  &{X_{n+1}} \ar@{-->}[r]^-{} &}
$$
Using similar arguments as in the proof of Lemma \ref{lem1}, we conclude
that the at least one of the following two distinguished $n$-exangles is non-split
$$\xymatrix{
U \ar[r]^{\delta_0}&L_1 \ar[r] & L_2 \ar[r]  & \cdots \ar[r] & {L_{n-1}}\ar[r]^{}  &{L_{n}}\ar[r]^{}  &{X_{n+1}} \ar@{-->}[r]^-{} &}$$
$$\xymatrix{
{V}\ar[r]^{\eta_0} &{N_1}\ar[r]^{}&{N_2} \ar[r]^{} &\cdot\cdot\cdot \ar[r] &{N_{n-1}}\ar[r]^{}  &{N_{n}}\ar[r]^{}  &{X_{n+1}} \ar@{-->}[r]^-{} &.
}$$
Without loss of generality, we assume that
$$\xymatrix{
U \ar[r]^{\delta_0}&L_1 \ar[r] & L_2 \ar[r]  & \cdots \ar[r] & {L_{n-1}}\ar[r]^{}  &{L_{n}}\ar[r]^{}  &{X_{n+1}} \ar@{-->}[r]^-{} &}$$ is non-split.
By Lemma \ref{ml}, there is a non-split distinguished $n$-exangle
$$\xymatrix{
C_{\bullet}:~U \ar[r]^{\quad\lambda_0}&C_1 \ar[r]^{\lambda_1} & C_2 \ar[r]^{\lambda_2}  & \cdots \ar[r]^{\lambda_{n-2}} & {C_{n-1}}\ar[r]^{\lambda_{n-1}}  &{C_{n}}\ar[r]^{\lambda_{n}}  &{X_{n+1}} \ar@{-->}[r]^-{} &}$$\
with $\lambda_1,\lambda_2,\cdots,\lambda_{n-1}$ in ${\rm rad_{\C}}$.
By {\rm (R0)} and the dual of Lemma \ref{a2}, we have the following commutative diagram
$$\xymatrix{
X_0\ar[r]^{\alpha_0}\ar@{}[dr] \ar@{-->}[d]^{} &X_1 \ar[r]^{\alpha_1} \ar@{}[dr]\ar@{-->}[d]^{}&X_2 \ar[r]^{\alpha_2} \ar@{}[dr]\ar@{-->}[d]^{}&\cdot\cdot\cdot \ar[r]^{\alpha_{n-2}} \ar@{}[dr] &X_{n-1} \ar[r]^{\alpha_{n-1}} \ar@{}[dr]\ar@{-->}[d]^{}& X_n\ar[r]^{\alpha_{n}} \ar@{}[dr]\ar[d]^{\binom{0}{1}}&{X_{n+1}} \ar@{}[dr]\ar@{=}[d]^{} \ar@{-->}[r]^-{} &\\
U\oplus V\ar[r]^{\;\;(u,~v)}\ar@{}[dr] \ar[d]^{p_1} &M_1 \ar[r]^{} \ar@{}[dr]\ar[d]^{\varphi_1}&M_2 \ar[r]^{} \ar@{}[dr]\ar[d]^{}&\cdot\cdot\cdot \ar[r]\ar@{}[dr] &M_{n-1} \ar[r]^{} \ar@{}[dr]\ar[d]^{}&U_n\oplus X_n\ar[r]^{(\beta_n,\alpha_n)} \ar@{}[dr]\ar[d]^{}&{X_{n+1}} \ar@{}[dr]\ar@{=}[d]^{} \ar@{-->}[r]^-{} &\\
U\ar[r]^{\delta_0}\ar@{}[dr] \ar@{=}[d]^{} &L_1 \ar[r]^{} \ar@{}[dr]\ar@{-->}[d]^{}&L_2 \ar[r]^{} \ar@{}[dr]\ar@{-->}[d]^{}&\cdot\cdot\cdot \ar[r]\ar@{}[dr] &L_{n-1} \ar[r]^{} \ar@{}[dr]\ar@{-->}[d]^{}&L_n\ar[r]^{} \ar@{}[dr]\ar@{-->}[d]^{}&{X_{n+1}} \ar@{}[dr]\ar@{=}[d]^{} \ar@{-->}[r]^-{} &\\
U \ar[r]^{\quad\lambda_0}&C_1 \ar[r]^{\lambda_1} & C_2 \ar[r]^{\lambda_2}  & \cdots \ar[r]^{\lambda_{n-2}} & {C_{n-1}}\ar[r]^{\lambda_{n-1}}  &{C_{n}}\ar[r]^{\lambda_{n}}  &{X_{n+1}} \ar@{-->}[r]^-{} &}
$$
of distinguished $n$-exangles. This shows that $X_{\bullet}>C_{\bullet}$.

By {\rm (R0)} and the dual of Lemma \ref{a2}, we have the following commutative diagram
$$\xymatrix{
U_0\ar[r]^{\beta_0}\ar@{}[dr] \ar@{-->}[d]^{} &U_1 \ar[r]^{\beta_1} \ar@{}[dr]\ar@{-->}[d]^{}&U_2 \ar[r]^{\beta_2} \ar@{}[dr]\ar@{-->}[d]^{}&\cdot\cdot\cdot \ar[r]^{\beta_{n-2}} \ar@{}[dr] &U_{n-1} \ar[r]^{\beta_{n-1}} \ar@{}[dr]\ar@{-->}[d]^{}& U_n\ar[r]^{\beta_{n}} \ar@{}[dr]\ar[d]^{\binom{1}{0}}&{X_{n+1}} \ar@{}[dr]\ar@{=}[d]^{} \ar@{-->}[r]^-{} &\\
U\oplus V\ar[r]^{\;\;(u,~v)}\ar@{}[dr] \ar[d]^{p_1} &M_1 \ar[r]^{} \ar@{}[dr]\ar[d]^{\varphi_1}&M_2 \ar[r]^{} \ar@{}[dr]\ar[d]^{}&\cdot\cdot\cdot \ar[r]\ar@{}[dr] &M_{n-1} \ar[r]^{} \ar@{}[dr]\ar[d]^{}&U_n\oplus X_n\ar[r]^{(\beta_n,\alpha_n)} \ar@{}[dr]\ar[d]^{}&{X_{n+1}} \ar@{}[dr]\ar@{=}[d]^{} \ar@{-->}[r]^-{} &\\
U\ar[r]^{\delta_0}\ar@{}[dr] \ar@{=}[d]^{} &L_1 \ar[r]^{} \ar@{}[dr]\ar@{-->}[d]^{}&L_2 \ar[r]^{} \ar@{}[dr]\ar@{-->}[d]^{}&\cdot\cdot\cdot \ar[r]\ar@{}[dr] &L_{n-1} \ar[r]^{} \ar@{}[dr]\ar@{-->}[d]^{}&L_n\ar[r]^{} \ar@{}[dr]\ar@{-->}[d]^{}&{X_{n+1}} \ar@{}[dr]\ar@{=}[d]^{} \ar@{-->}[r]^-{} &\\
U \ar[r]^{\quad\lambda_0}&C_1 \ar[r]^{\lambda_1} & C_2 \ar[r]^{\lambda_2}  & \cdots \ar[r]^{\lambda_{n-2}} & {C_{n-1}}\ar[r]^{\lambda_{n-1}}  &{C_{n}}\ar[r]^{\lambda_{n}}  &{X_{n+1}} \ar@{-->}[r]^-{} &}
$$
of distinguished $n$-exangles. This shows that $U_{\bullet}>C_{\bullet}$.\qed
\medskip

\begin{lemma}\label{y8}
Let $\C$ be a locally finite $n$-exangulated category.
\begin{itemize}
\item[\rm (1)] If $X_{n+1}\in\ind(\C)$ is a non-projective object, then $S(X_{n+1})$ has a minimal element.
\item[ \rm (2)] If $Y_{0}\in\ind(\C)$ is a non-injective object, then $T(Y_{0})$ has a minimal element.

\end{itemize}
\end{lemma}\proof We just prove the first statement, the second statement proves similarly.

Since $X_{n+1}\in\ind(\C)$ is a non-projective, there exists an object $X_{0}\in\C$, such that $\E(X_{n+1},X_{0})\neq 0$ by Lemma \ref{def2}. That is to say, there is a non-split distinguished $n$-exangle:$$X_{\bullet}:X_0\xrightarrow{\alpha_0}X_1\xrightarrow{\alpha_1}X_2\xrightarrow{\alpha_2}\cdots\xrightarrow{\alpha_{n-1}}
X_n\xrightarrow{\alpha_n}X_{n+1}\overset{\delta}{\dashrightarrow}.$$
Since $\C$ is a Krull-Schmidt category, we decompose $X_n$ into a direct sum of indecomposable objects $X_n=\bigoplus\limits_{k=1}^rB_k$.
Thus $X_{\bullet}$ can be written as
$$X_{\bullet}:X_0\xrightarrow{\alpha_0}X_1\xrightarrow{\alpha_1}X_2\xrightarrow{\alpha_2}\cdots\xrightarrow{\alpha_{n-1}}
\bigoplus\limits_{k=1}^rB_k\xrightarrow{(b_1,b_2,\cdots,b_r)}X_{n+1}\overset{}{\dashrightarrow}.$$

where $b_k\in{\rm rad}_{\C}(B_k,X_{n+1})$, $k=1,2,\cdots,r$.

Since $\C$ is locally finite, there are only finite many objects $X_i\in\ind(\C),~i=1,2,\cdots,m$
such that $\Hom_{\C}(X_i,X_{n+1})\neq 0$.
We assume that $\lambda_{ij},~1\leq j\leq q_i$ form a basis of the $k$-vector space ${\rm rad}_{\C}(B_k,X_{n+1})$.
Put $M:=(\bigoplus\limits_{k=1}^rB_k)\oplus (\bigoplus\limits_{i=1}^m(X_i)^{\oplus q_i})$, we consider the morphism
$$\gamma:=(b_1,b_2,\cdots,b_r,\lambda_{11},\cdots,\lambda_{ij},\cdots,\lambda_{mq_m})\in{\rm rad}_{\C}(M,X_{n+1})$$
which is not split epimorphism. By {\rm (EA2)}, we deduce that there is a distinguished $n$-exangle in $\C$ as follows:
$$M_{\bullet}:~B\xrightarrow{} M_1\xrightarrow{}
 M_2\xrightarrow{}\cdots\xrightarrow{}M_{n-1}
                                     \xrightarrow{}M\xrightarrow{\gamma}X_{n+1}\overset{}{\dashrightarrow}.$$

Thus $M_{\bullet}$ is non-split since $\gamma$ is not split epimorphism.
Without loss of generality, we can assume that $B=U\oplus V$ where $U$ and $V$ are indecomposable.
For the morphism $p_1=(1,0)\colon U\oplus V\to U$, by  (EA2$\op$), we can observe that $(p_1,\id_{X_{n+1}})$ has a {\it good lift} $g^{\mr}=(p_1,\varphi_1,\ldots,\varphi_n,\id_{X_{n+1}})$, that is,
there exists the following commutative diagram of distinguished $n$-exangles
$$\xymatrix{
U\oplus V\ar[r]^{\;\;(u,~v)}\ar@{}[dr] \ar[d]^{p_1} &M_1 \ar[r]^{} \ar@{}[dr]\ar@{-->}[d]^{\varphi_1}&M_2 \ar[r]^{} \ar@{}[dr]\ar@{-->}[d]^{}&\cdot\cdot\cdot \ar[r]\ar@{}[dr] &M_{n-1} \ar[r]^{} \ar@{}[dr]\ar@{-->}[d]^{}&{M}\ar[r]^{\gamma} \ar@{}[dr]\ar@{-->}[d]^{}&{X_{n+1}} \ar@{}[dr]\ar@{=}[d]^{} \ar@{-->}[r]^-{} &\\
{U}\ar[r]^{\theta_0} &{L_1}\ar[r]^{}&{L_2} \ar[r]^{} &\cdot\cdot\cdot \ar[r] &{L_{n-1}}\ar[r]^{}  &{L_{n}}\ar[r]^{}  &{X_{n+1}} \ar@{-->}[r]^-{} &}
$$

Similarly, for the morphism $p_2=(0,1)\colon U\oplus V\to V$, there exists
the following commutative diagram of distinguished $n$-exangles
$$\xymatrix{
U\oplus V\ar[r]^{\;\;(u,~v)}\ar@{}[dr] \ar[d]^{p_2} &M_1 \ar[r]^{} \ar@{}[dr]\ar@{-->}[d]^{m_1}&M_2 \ar[r]^{} \ar@{}[dr]\ar@{-->}[d]^{}&\cdot\cdot\cdot \ar[r]\ar@{}[dr] &M_{n-1} \ar[r]^{} \ar@{}[dr]\ar@{-->}[d]^{}&{M}\ar[r]^{\gamma} \ar@{}[dr]\ar@{-->}[d]^{}&{X_{n+1}} \ar@{}[dr]\ar@{=}[d]^{} \ar@{-->}[r]^-{} &\\
{V}\ar[r]^{\eta_0} &{N_1}\ar[r]^{}&{N_2} \ar[r]^{} &\cdot\cdot\cdot \ar[r] &{N_{n-1}}\ar[r]^{}  &{N_{n}}\ar[r]^{}  &{X_{n+1}} \ar@{-->}[r]^-{} &}
$$
Using similar arguments as in the proof of Lemma \ref{lem1}, we conclude
that the at least one of the following two distinguished $n$-exangles is non-split
$$\xymatrix{
U \ar[r]^{\theta_0}&L_1 \ar[r] & L_2 \ar[r]  & \cdots \ar[r] & {L_{n-1}}\ar[r]^{}  &{L_{n}}\ar[r]^{}  &{X_{n+1}} \ar@{-->}[r]^-{} &}$$
$$\xymatrix{
{V}\ar[r]^{\eta_0} &{N_1}\ar[r]^{}&{N_2} \ar[r]^{} &\cdot\cdot\cdot \ar[r] &{N_{n-1}}\ar[r]^{}  &{N_{n}}\ar[r]^{}  &{X_{n+1}} \ar@{-->}[r]^-{} &.
}$$
Without loss of generality, we assume that
$$\xymatrix{
U \ar[r]^{\theta_0}&L_1 \ar[r] & L_2 \ar[r]  & \cdots \ar[r] & {L_{n-1}}\ar[r]^{}  &{L_{n}}\ar[r]^{}  &{X_{n+1}} \ar@{-->}[r]^-{} &}$$ is non-split.
By Lemma \ref{ml}, we can find a non-split distinguished $n$-exangle
$$\xymatrix{
C_{\bullet}:~U \ar[r]^{\quad\omega_0}&C_1 \ar[r]^{\omega_1} & C_2 \ar[r]^{\omega_2}  & \cdots \ar[r]^{\omega_{n-2}} & {C_{n-1}}\ar[r]^{\omega_{n-1}}  &{C_{n}}\ar[r]^{\omega_{n}}  &{X_{n+1}} \ar@{-->}[r]^-{} &}$$\
with $\omega_1,\omega_2,\cdots,\omega_{n-1}$ in ${\rm rad_{\C}}$. Then $C_{\bullet}\in S(X_{n+1})$.
By (R0), we have the following commutative diagram
$$\xymatrix{
U\oplus V\ar[r]^{\;\;(u,~v)}\ar@{}[dr] \ar[d]^{p_1} &M_1 \ar[r]^{} \ar@{}[dr]\ar[d]^{\varphi_1}&M_2 \ar[r]^{} \ar@{}[dr]\ar[d]^{}&\cdot\cdot\cdot \ar[r]\ar@{}[dr] &M_{n-1} \ar[r]^{} \ar@{}[dr]\ar[d]^{}&{M}\ar[r]^{\gamma}\ar@{}[dr]\ar[d]^{}&{X_{n+1}} \ar@{}[dr]\ar@{=}[d]^{} \ar@{-->}[r]^-{} &\\U\ar[r]^{\theta_0}\ar@{}[dr] \ar@{=}[d]^{} &L_1 \ar[r]^{} \ar@{}[dr]\ar@{-->}[d]^{}&L_2 \ar[r]^{} \ar@{}[dr]\ar@{-->}[d]^{}&\cdot\cdot\cdot \ar[r]\ar@{}[dr] &L_{n-1} \ar[r]^{} \ar@{}[dr]\ar@{-->}[d]^{}&L_n\ar[r]^{} \ar@{}[dr]\ar@{-->}[d]^{}&{X_{n+1}} \ar@{}[dr]\ar@{=}[d]^{} \ar@{-->}[r]^-{} &\\
U \ar[r]^{\quad\omega_0}&C_1 \ar[r]^{\omega_1} & C_2 \ar[r]^{\omega_2}  & \cdots \ar[r]^{\omega_{n-2}} & {C_{n-1}}\ar[r]^{\omega_{n-1}}  &{C_{n}}\ar[r]^{\omega_{n}}  &{X_{n+1}} \ar@{-->}[r]^-{} &}
$$
of distinguished $n$-exangles.
For any $D_{\bullet}\in S(X_{n+1})$, it can be written as
$$\xymatrix {D_{\bullet}:~D\xrightarrow{~}D_1 \xrightarrow{~} D_2 \xrightarrow{~} \cdots
  \xrightarrow{~} D_{n-1} \xrightarrow{~} \bigoplus\limits_{i=1}^pH_i \xrightarrow{d=(d_1,d_2,\cdots,d_p)}X_{n+1}\overset{}{\dashrightarrow}}$$
with $H_i\in\ind(\C), d_i\in{\rm rad}_{\C}(H_i,X_{n+1})$, $i=1,2,\cdots,p$.
Since $D_{\bullet}\in S(X_{n+1})$ is non-split, $d$ is not split epimorphism implies that
$d\in{\rm rad}_{\C}(\bigoplus\limits_{i=1}^pH_i ,X_{n+1})$.
By the definitions of $\lambda_{ij}$ and $\gamma$, there exists a morphism
$\rho\colon \bigoplus\limits_{i=1}^pH_i\to M$ such that $d=\gamma\rho$.
By the dual of Lemma \ref{a2}, we have the following commutative diagram
\vspace{-2mm}
$$\xymatrix{
D\ar[r]\ar@{-->}[d]& D_1 \ar[r]\ar@{-->}[d] & D_2 \ar[r]\ar@{-->}[d]& \cdots \ar[r]&\ar[r]\ar@{-->}[d]D_{d-1}&\bigoplus\limits_{i=1}^pH_i \ar[r]^{\quad d}\ar[d]^{\rho}&X_{n+1}\ar@{}[dr]\ar@{=}[d]^{} \ar@{-->}[r]^-{} &\\
B\ar[r]&M_1 \ar[r] & M_2 \ar[r]  & \cdots \ar[r] &M_{d-1}\ar[r]& M \ar[r]^{\gamma}&X_{n+1}\ar@{-->}[r]^-{} &}$$
of distinguished $n$-exangles, where $B=U\oplus V$.
Thus we get the following commutative diagram
$$\xymatrix{
D\ar[r]\ar@{-->}[d]& D_1 \ar[r]\ar@{-->}[d] & D_2 \ar[r]\ar@{-->}[d]& \cdots \ar[r]&\ar[r]\ar@{-->}[d]D_{d-1}&\bigoplus\limits_{i=1}^pH_i \ar[r]^{\quad d}\ar[d]^{}&X_{n+1}\ar@{}[dr]\ar@{=}[d]^{} \ar@{-->}[r]^-{} &\\
U \ar[r]^{\quad\omega_0}&C_1 \ar[r]^{\omega_1} & C_2 \ar[r]^{\omega_2}  & \cdots \ar[r]^{\omega_{n-2}} & {C_{n-1}}\ar[r]^{\omega_{n-1}}  &{C_{n}}\ar[r]^{\omega_{n}}  &{X_{n+1}} \ar@{-->}[r]^-{} &}$$
of distinguished $n$-exangles. This shows that $C_{\bullet}$ is a minimal element in $S(X)$.  \qed

We are now ready to state and prove our main result.

\begin{theorem}\label{main}
Let $\C$ be a locally finite $n$-exangulated category. If $X_{n+1}\in\ind(\C)$ is a non-projective object, then there exists an Auslander-Reiten $n$-exangle ending at $X_{n+1}$, and if $Y_{0}\in\ind(\C)$ is a non-injective object, then there exists an Auslander-Reiten $n$-exangle starting at $Y_{0}$.
In this case, we say that $\C$ has Auslander-Reiten $n$-exangles.
\end{theorem}
\proof Since $X_{n+1}\in\ind(\C)$, by Lemma \ref{lem1} we know that the set $S(X_{n+1})$ is non-empty.
Thus by Lemma \ref{y8}, there is a  distinguished $n$-exangle
$$X_{\bullet}:X_0\xrightarrow{\alpha_0}X_1\xrightarrow{\alpha_1}X_2\xrightarrow{\alpha_2}\cdots\xrightarrow{\alpha_{n-1}}
X_n\xrightarrow{\alpha_n}X_{n+1}\overset{\delta}{\dashrightarrow}$$
where $\alpha_1,\alpha_2,\cdots,\alpha_{{n-1}}\in {\rm rad}_{\C}$ and $X_0\in\ind(\C)$, such that $X_{\bullet}$ is a minimal element in $S(X_{n+1})$. Then $\End(X_0)$ is local.

We want to prove that $X_{\bullet}$ is an Auslander-Reiten $n$-exangle, by Lemma \ref{r1},
it is enough to show that $\alpha_n$ is right almost split.

Assume that $g\colon M_{n+1}\to X_{n+1}$ is not a split epimorphism, we claim that $g$ factors through $\alpha_{{n}}$.
By {\rm (EA2)}, we can observe that $(\id_{X_{0}},g)$ has a {\it good lift} $g^{\mr}=(\id_{X_{0}},\varphi_1,\ldots,\varphi_n,g)$, that is,
there exists the following commutative diagram of distinguished $n$-exangles
$$\xymatrix{
X_0\ar[r]^{\gamma_0}\ar@{}[dr] \ar@{=}[d]^{} &B_1 \ar[r]^{\gamma_1} \ar@{}[dr]\ar@{-->}[d]^{\varphi_1}&B_2 \ar[r]^{\gamma_2} \ar@{}[dr]\ar@{-->}[d]^{\varphi_2}&\cdot\cdot\cdot \ar[r]\ar@{}[dr] &B_n \ar[r]^{\gamma_n} \ar@{}[dr]\ar@{-->}[d]^{\varphi_n}&{M_{n+1}} \ar@{}[dr]\ar[d]^{g} \ar@{-->}[r]^-{g^{\ast}\delta} &\\
{X_0}\ar[r]^{\alpha_0} &{X_1}\ar[r]^{\alpha_1}&{X_2} \ar[r]^{\alpha_2} &\cdot\cdot\cdot \ar[r] &{X _n}\ar[r]^{\alpha_n}  &{X_{n+1}} \ar@{-->}[r]^-{\delta} &}
$$
such that
$$N_{\bullet}:~B_1\xrightarrow{} N_1\xrightarrow{}
 N_2\xrightarrow{}\cdots\xrightarrow{}N_{n-1}
                                     \xrightarrow{}{M_{n+1}} \oplus X_n\xrightarrow{\left(
                                       \begin{smallmatrix}
                                         g,&\alpha_{n}
                                       \end{smallmatrix}
                                     \right)}X_{n+1}\overset{(\gamma_0)_{\ast}\delta}{\dashrightarrow}$$
 is a distinguished $n$-exangle in $\C$, where $N_i=B_{i+1}\oplus X_{i},~i=1,2,\cdots,n-1$.
Since $g$ and $\alpha_n$ are not split epimorphisms, we have that
$(g,\alpha_n)$ is also not split epimorphism by using similar arguments as in the proof of Lemma \ref{y7}. That is, $N_{\bullet}$ is non-split.

Without loss of generality, we can assume that $B_1=U\oplus V$  where $U$ and $V$ are indecomposable. For the morphism $p_1=(1,0)\colon U\oplus V\to U$, by  (EA2$\op$), we can observe that $(p_1,\id_{X_{n+1}})$ has a {\it good lift} $h^{\mr}=(p_1,\phi_1,\ldots,\phi_n,\id_{X_{n+1}})$, that is,
there exists the following commutative diagram of distinguished $n$-exangles
$$\xymatrix{
U\oplus V\ar[r]^{\;\;(u,~v)}\ar@{}[dr] \ar[d]^{p_1} &N_1 \ar[r]^{} \ar@{}[dr]\ar@{-->}[d]^{\phi_1}&N_2 \ar[r]^{} \ar@{}[dr]\ar@{-->}[d]^{}&\cdot\cdot\cdot \ar[r]\ar@{}[dr] &N_{n-1} \ar[r]^{} \ar@{}[dr]\ar@{-->}[d]^{}&{M_{n+1}} \oplus X_n\ar[r]^{} \ar@{}[dr]\ar@{-->}[d]^{}&{X_{n+1}} \ar@{}[dr]\ar@{=}[d]^{} \ar@{-->}[r]^-{} &\\
{U}\ar[r]^{\delta_0} &{L_1}\ar[r]^{}&{L_2} \ar[r]^{} &\cdot\cdot\cdot \ar[r] &{L_{n-1}}\ar[r]^{}  &{L_{n}}\ar[r]^{}  &{X_{n+1}} \ar@{-->}[r]^-{} &}
$$

Similarly, for the morphism $p_2=(0,1)\colon U\oplus V\to V$, there exists
the following commutative diagram of distinguished $n$-exangles
$$\xymatrix{
U\oplus V\ar[r]^{\;\;(u,~v)}\ar@{}[dr] \ar[d]^{p_2} &N_1 \ar[r]^{} \ar@{}[dr]\ar@{-->}[d]^{q_1}&N_2 \ar[r]^{} \ar@{}[dr]\ar@{-->}[d]^{}&\cdot\cdot\cdot \ar[r]\ar@{}[dr] &N_{n-1} \ar[r]^{} \ar@{}[dr]\ar@{-->}[d]^{}&{M_{n+1}} \oplus X_n\ar[r]^{} \ar@{}[dr]\ar@{-->}[d]^{}&{X_{n+1}} \ar@{}[dr]\ar@{=}[d]^{} \ar@{-->}[r]^-{} &\\
{V}\ar[r]^{\eta_0} &{Q_1}\ar[r]^{}&{Q_2} \ar[r]^{} &\cdot\cdot\cdot \ar[r] &{Q_{n-1}}\ar[r]^{}  &{Q_{n}}\ar[r]^{}  &{X_{n+1}} \ar@{-->}[r]^-{} &}
$$
Using similar arguments as in the proof of Lemma \ref{lem1}, we conclude
that the at least one of the following two distinguished $n$-exangles is non-split
$$\xymatrix{
{U}\ar[r]^{\delta_0} &{L_1}\ar[r]^{}&{L_2} \ar[r]^{} &\cdot\cdot\cdot \ar[r] &{L_{n-1}}\ar[r]^{}  &{L_{n}}\ar[r]^{}  &{X_{n+1}} \ar@{-->}[r]^-{} &}$$
$$\xymatrix{
{V}\ar[r]^{\eta_0} &{Q_1}\ar[r]^{}&{Q_2} \ar[r]^{} &\cdot\cdot\cdot \ar[r] &{Q_{n-1}}\ar[r]^{}  &{Q_{n}}\ar[r]^{}  &{X_{n+1}} \ar@{-->}[r]^-{} &.
}$$
Without loss of generality, we assume that
$$\xymatrix{
{U}\ar[r]^{\delta_0} &{L_1}\ar[r]^{}&{L_2} \ar[r]^{} &\cdot\cdot\cdot \ar[r] &{L_{n-1}}\ar[r]^{}  &{L_{n}}\ar[r]^{}  &{X_{n+1}} \ar@{-->}[r]^-{} &}$$ is non-split.
By Lemma \ref{ml}, we can find a non-split distinguished $n$-exangle
$$\xymatrix{
C_{\bullet}:~U \ar[r]^{\quad\lambda_0}&C_1 \ar[r]^{\lambda_1} & C_2 \ar[r]^{\lambda_2}  & \cdots \ar[r]^{\lambda_{n-2}} & {C_{n-1}}\ar[r]^{\lambda_{n-1}}  &{C_{n}}\ar[r]^{\lambda_{n}}  &{X_{n+1}} \ar@{-->}[r]^-{} &}$$\
with $\lambda_1,\lambda_2,\cdots,\lambda_{n-1}$ in ${\rm rad_{\C}}$.
By {\rm (R0)} and the dual of Lemma \ref{a2}, we have the following commutative diagram
$$\xymatrix{
X_0\ar[r]^{\alpha_0}\ar@{}[dr] \ar@{-->}[d]^{} &X_1 \ar[r]^{\alpha_1} \ar@{}[dr]\ar@{-->}[d]^{}&X_2 \ar[r]^{\alpha_2} \ar@{}[dr]\ar@{-->}[d]^{}&\cdot\cdot\cdot \ar[r]^{\alpha_{n-2}} \ar@{}[dr] &X_{n-1} \ar[r]^{\alpha_{n-1}} \ar@{}[dr]\ar@{-->}[d]^{}& X_n\ar[r]^{\alpha_{n}} \ar@{}[dr]\ar[d]^{\binom{0}{1}}&{X_{n+1}} \ar@{}[dr]\ar@{=}[d]^{} \ar@{-->}[r]^-{} &\\U\oplus V\ar[r]^{\;\;(u,~v)}\ar@{}[dr] \ar[d]^{p_1} &N_1 \ar[r]^{} \ar@{}[dr]\ar[d]^{\phi_1}&N_2 \ar[r]^{} \ar@{}[dr]\ar[d]^{}&\cdot\cdot\cdot \ar[r]\ar@{}[dr] &N_{n-1} \ar[r]^{} \ar@{}[dr]\ar[d]^{}&{M_{n+1}} \oplus X_n\ar[r]^{(g,\alpha_n)} \ar@{}[dr]\ar[d]^{}&{X_{n+1}} \ar@{}[dr]\ar@{=}[d]^{} \ar@{-->}[r]^-{} &\\U\ar[r]^{\delta_0}\ar@{}[dr] \ar@{=}[d]^{} &L_1 \ar[r]^{} \ar@{}[dr]\ar@{-->}[d]^{}&L_2 \ar[r]^{} \ar@{}[dr]\ar@{-->}[d]^{}&\cdot\cdot\cdot \ar[r]\ar@{}[dr] &L_{n-1} \ar[r]^{} \ar@{}[dr]\ar@{-->}[d]^{}&L_n\ar[r]^{} \ar@{}[dr]\ar@{-->}[d]^{}&{X_{n+1}} \ar@{}[dr]\ar@{=}[d]^{} \ar@{-->}[r]^-{} &\\
U \ar[r]^{\quad\lambda_0}&C_1 \ar[r]^{\lambda_1} & C_2 \ar[r]^{\lambda_2}  & \cdots \ar[r]^{\lambda_{n-2}} & {C_{n-1}}\ar[r]^{\lambda_{n-1}}  &{C_{n}}\ar[r]^{\lambda_{n}}  &{X_{n+1}} \ar@{-->}[r]^-{} &}
$$
of distinguished $n$-exangles.
We obtain that $X_{\bullet}>C_{\bullet}$ implies that $X_{\bullet}\sim C_{\bullet}$ since $X_{\bullet}$ is the minimal element in $S(X_{n+1})$. Thus there exists the following commutative diagram
$$\xymatrix{
U\ar[r]^{\lambda_0}\ar@{}[dr] \ar[d]^{} &C_1 \ar[r]^{\lambda_1} \ar@{}[dr]\ar[d]^{}&C_2 \ar[r]^{\lambda_2} \ar@{}[dr]\ar[d]^{}&\cdot\cdot\cdot \ar[r]^{\lambda_{n-2}}\ar@{}[dr] &C_{n-1} \ar[r]^{\lambda_{n-1}}\ar@{}[dr] \ar[d]^{} &C_n \ar[r]^{\lambda_n} \ar@{}[dr]\ar[d]^{}&X_{n+1} \ar@{}[dr]\ar@{=}[d]^{} \ar@{-->}[r]^-{} &\\
{X_0}\ar[r]^{\alpha_0} &{X_1}\ar[r]^{\alpha_1}&{X_2} \ar[r]^{\alpha_2} &\cdot\cdot\cdot \ar[r]^{\alpha_{n-2}} &{X_{n-1} }\ar[r]^{\alpha_{n-1} } &{X_n}\ar[r]^{\alpha_n}  &{X_{n+1}} \ar@{-->}[r]^-{} &}
$$
of distinguished $n$-exangles.
Hence we get the following commutative diagram
\vspace{-4mm}
$$\xymatrix{
U\oplus V\ar[r]^{\;\;(u,~v)}\ar@{}[dr] \ar[d]^{} &N_1 \ar[r]^{} \ar@{}[dr]\ar[d]^{} &N_2 \ar[r]^{} \ar@{}[dr]\ar[d]^{} &\cdot\cdot\cdot \ar[r]\ar@{}[dr] &N_{n-1} \ar[r]^{} \ar@{}[dr]\ar[d]^{} &{M_{n+1}} \oplus X_n\ar[r]^{\quad (g,\alpha_n)} \ar@{}[dr]\ar[d]^{(a,b)} &{X_{n+1}} \ar@{}[dr]\ar@{=}[d]^{} \ar@{-->}[r]^-{} &\\
{X_0}\ar[r]^{\alpha_0} &{X_1}\ar[r]^{\alpha_1}&{X_2} \ar[r]^{\alpha_2} &\cdot\cdot\cdot \ar[r]^{\alpha_{n-2}} &{X_{n-1} }\ar[r]^{\alpha_{n-1} } &{X_n}\ar[r]^{\alpha_n}  &{X_{n+1}} \ar@{-->}[r]^-{} &}
$$
of distinguished $n$-exangles. It follows that $g=\alpha_na$.
This shows that $\alpha_n$ is right almost split.

Similarly, we can show that if $Y_{0}\in\ind(\C)$, then there exists an Auslander-Reiten $n$-exangle starting at $Y_{0}$.
Thus $\C$ has Auslander-Reiten $n$-exangles. \qed
\medskip

By applying Theorem \ref{main} to $(n+2)$-angulated category, we have the following.

\begin{corollary}\rm\cite[Theorem 1.1]{Z2}
Let $\C$ be a locally finite $(n+2)$-angulated category. Then $\C$ has Auslander-Reiten $(n+2)$-angles.
\end{corollary}

By applying Theorem \ref{main} to $n$-abelian categories, we have the following.
\begin{corollary}\rm\cite[Theorem 1.1]{XLW}
Let $\C$ be a locally finite $n$-abelian category. Then $\C$ has $n$-Auslander-Reiten sequences.
\end{corollary}

By applying Theorem \ref{main} to $n$-exact categories, we have the following.

\begin{corollary}
Let $\C$ be a locally finite $n$-exact category. Then $\C$ has $n$-Auslander-Reiten sequences.
\end{corollary}

\begin{remark}
As a special case of Theorem \ref{main} when $n=1$, that is, if $\C$ is a locally finite
extriangulated category, then $\C$ has Auslander-Reiten  $\E$-triangles, see \cite{ZZ4}.
\end{remark}

\begin{remark}
If $\C$ is a locally finite
triangulated category, then $\C$ has Auslander-Reiten triangles, see \cite{XZ1,XZ2}.
\end{remark}

\textbf{Jian He}\\
Department of Mathematics, Nanjing University, 210093 Nanjing, Jiangsu, P. R. China\\
E-mail: \textsf{jianhe30@163.com}\\[1mm]
\textbf{Jiangsheng Hu}\\
School of Mathematics and Physics, Jiangsu University of Technology,
 Changzhou, Jiangsu 213001, P. R. China.\\
E-mail: \textsf{jiangshenghu@jsut.edu.cn}\\[1mm]
\textbf{Dongdong Zhang}\\
Department of Mathematics, Zhejiang Normal University,
321004 Jinhua, Zhejiang, P. R. China.\\
E-mail: \textsf{zdd@zjnu.cn}\\[1mm]
\textbf{Panyue Zhou}\\
College of Mathematics, Hunan Institute of Science and Technology, 414006 Yueyang, Hunan, P. R. China.\\
E-mail: \textsf{panyuezhou@163.com}

\end{document}